\documentclass[12pt]{article} 
\usepackage{latexsym}
\usepackage{amssymb}
\usepackage{euscript}

\let\cal\mathcal

\def\ie{{\it i.e. }}
\def\eg{{\it e.g. }}
\def\={\ = \ }
\def\la{\langle}
\def\ra{\rangle}

\def\be{\setcounter{equation}{\value{theorem}} \begin{equation}}
\def\ee{\end{equation} \addtocounter{theorem}{1}}
\def\bp{{\sc Proof: }}
\def\ep{{}{\hfill $\Box$} \vskip 5pt \par}
\def\bl{\begin{lemma}}
\def\el{\end{lemma}}
\def\bt{\begin{theorem}}
\def\et{\end{theorem}}
\def\bprop{\begin{prop}}
\def\eprop{\end{prop}}

\def\sz{Szeg\H{o} }

\def\h{{\cal H}}
\def\hk{{\cal H}_k}

\def\L{{\cal L}}

\def\D{{\mathbb D}}

\def\C{{\mathbb C}}

\def\={\ = \ }

\def\l{\lambda}
\def\G{\Gamma}

\def\V{\Vert}
\def\i{\infty}

\def\hd{H^\infty(\D)}
\def\hb{H^\infty(\D^2)}

\advance\topmargin by -.6in
\advance\textheight by 1.2in
\advance\oddsidemargin by -.5in
\advance\textwidth by 1.05in

\newtheorem{theorem}{Theorem} [section]
\newtheorem{prop}[theorem]{Proposition}
\newtheorem{lemma}[theorem]{Lemma}

\newtheorem{question}[theorem]{Question}

\begin{document}
\setlength{\baselineskip}{21pt}
\title{ Interpolating Sequences on the Bidisk}
\author{Jim Agler
\thanks{Partially supported by the National Science Foundation}\\
U.C. San Diego\\
La Jolla, California 92093
\and
John E. M\raise.5ex\hbox{c}Carthy\\
Washington University\\
St. Louis, Missouri 63130}
\date{}

\bibliographystyle{plain}

\maketitle


\baselineskip = 18pt

\setcounter{section}{-1}
\section{Introduction}
By an {\it interpolating sequence} for $\hd$, the algebra of bounded analytic
functions on the unit disk $\D$, we mean a sequence $\{ \l_i \}_{i=1}^\i$ in
$\D$ with the property that, for any bounded sequence of complex numbers
$\{ w_i \}_{i=1}^\i$, there is a function $f$  in $\hd$ such that $f(\l_i) =
w_i$ for each $i$. These sequences were characterized by L.~Carleson in 1958
\cite{car58} (Theorem~\ref{car} below).
To state Carleson's theorem, we need a few definitions.

{\bf Definition:} The sequence $\{ \l_i \}_{i=1}^\i $ is {\it weakly
separated} if there exists a constant $M$ such that, whenever
$i \neq j$, there exists a function $\phi$ in $\hd$ of norm at most
$M$ that satisfies $\phi(\l_i) = 1$ and $\phi(\l_j) = 0$.
The sequence is {\it strongly separated} if, for each $i$ there is
a $\phi$ in $\hd$ of norm at most
$M$ that satisfies $\phi(\l_i) = 1$ and $\phi(\l_j) = 0$
for all $j \neq i$.

When dealing with $\hd$, it is natural to consider it as the multiplier
algebra of the Hardy space $H^2$, the Hilbert space of holomorphic functions
on the disk with norm
$$
\| f \|^2_{H^2} \ := \
\frac{1}{2\pi}\, \lim_{r \nearrow 1} \int |f(re^{i\theta}) |^2 d\theta .
$$
The reproducing kernel for $H^2$ is the \sz kernel
$$
s(\l,z) \= \frac{1}{1-\bar \l z} .
$$
With respect to this kernel, the 
{\it associated Gram matrix} of the sequence
$\{ \l_i \}_{i=1}^\i $  is the Grammian of the normalized kernel
functions, \ie the matrix with $(i,j)$ entry given by
$$
G_{ij} \= \frac{s(\l_i,\l_j)}{\sqrt{ s(\l_i,\l_i) s(\l_j,\l_j)}} .
$$

We can now state Carleson's theorem \cite{car58}.

\bt {\bf [Carleson]}
\label{car}
Let $\{ \l_i \}_{i=1}^\i $ be a sequence in $\D$. Then the following are
equivalent:

(1) $\{ \l_i \}_{i=1}^\i $ is an interpolating sequence for $\hd$.

(2) The sequence is weakly separated and the associated Gram matrix is
bounded.

(3) The sequence is strongly separated.
\et

Condition (2) is normally written differently. Rather than saying the
associated Gram matrix is bounded,
one says instead that the measure
$ \displaystyle
\sum (1 - |\l_i|^2) \delta_{\l_i}$
is a Carleson measure, \ie there exists some constant $C$ such that, for
all $f$ in $H^2$, the inequality
$$
\sum_{i=1}^\i (1 - |\l_i|^2) | f(\l_i) |^2 \ \leq \
C^2\, \V f \V_{H^2}^2
$$
holds.
The fact that boundedness of the Gram matrix is equivalent to a particular
atomic measure being Carleson is true in great generality, and is proved in
Proposition~\ref{bgcm} below.

In 1987, B. Berndtsson, S.-Y. Chang and K.-C. Lin 
studied interpolating sequences for $H^\i(\D^d)$, the bounded analytic
functions on the polydisk \cite{bcl}. To state their result, let us first define
the Gleason distance between two points by
$$
\rho(\zeta,\l) \ := \
\sup \{|f(\zeta)| : \| f \|_{H^\i(\D^d)} \leq 1, f(\l) = 0 \} .
$$
Notice that the elementary theory of Blaschke products on the disk shows
that condition (3) in Theorem~\ref{car} is equivalent to:

{\it
(3)': There exists $c > 0$ such that 
$$
\prod_{j\neq i} \rho(\l_i, \l_j) \geq c
$$
for all $i$.
}

With respect to Lebesgue measure $\sigma$ on the distinguished boundary of
the polydisk, one can define a space $H^2(\sigma)$ as the closure of the
polynomials in $L^2(\sigma)$. This space has a reproducing kernel
$$
k(\l, z) \= \prod_{n=1}^d \frac{1}{1- \bar \l^n z^n},
$$
where we use $\l^n$ to denote the $n^{th}$ component of the point $\l$.

\bt
{\bf [Berndtsson,  Chang and Lin]}
Consider the three statements

(1) There exists $c > 0$ such that
$$
\prod_{j\neq i} \rho(\l_i, \l_j) \geq c
$$
for all $i$.

(2) The sequence $\{ \l_i \}_{i=1}^\i$ is an interpolating sequence for 
$H^\i(\D^d)$.

(3) The sequence $\{ \l_i \}_{i=1}^\i$ is weakly separated and the
associated Grammian with respect to Lebesgue measure $\sigma$ is bounded.

Then (1) implies (2) and (2) implies (3). Moreover the converse of both
these implications is false.
\et

For additional sufficient conditions that
guarantee that a sequence is interpolating, see the papers of E. Kronstadt
\cite{kro} and L. Huang \cite{hua}.

It is the purpose of this paper to give a characterization of interpolating
sequences for $\hb$.
To this end, we observe that whereas all pure cyclic isometries are unitarily
equivalent, so that in one variable one need only consider the \sz
kernel,
there are many non-equivalent pairs of pure cyclic isometries.
Thus when considering $\hb$ it is essential to consider a whole family of
kernels simultaneously.

Given a kernel $k$ on the bidisk, say it is {\it admissible}
 if
$$
(1 - \overline{\l^1} z^1) k(\l,z) \ \geq \ 0
$$
and
$$
(1 - \overline{\l^2} z^2) k(\l,z) \ \geq \ 0 ,
$$
\ie if
multiplication by each coordinate function is a contraction
(see Section~\ref{sec:ker} below).

For a given sequence $\{ \l_i \}_{i=1}^\i $,
the normalized Grammian of $k$,
$G^k$, is then the infinite matrix

$$
G^k_{ij} \=
\frac{k(\l_i,\l_j)}{\sqrt{ k(\l_i,\l_i) k(\l_j,\l_j)}}.
$$

Let $I$ denote the identity matrix $\delta_{ij}$,
and $J$ the matrix all of whose entries are $1$. 
We can now state our main result.

\bt
\label{thm:is}
Let $\{ \l_i \}_{i=1}^\i $ be a sequence in $\D^2$. The following are
equivalent:

(i) $\{ \l_i \}_{i=1}^\i $ is an interpolating sequence for $\hb$.

(ii) The following two conditions hold.

\noindent
(a)  For all admissible kernels $k$,
their normalized Grammians are uniformly bounded:
$$
 G^k  \ \leq \ M I 
$$
for some $M$.

(b) For all admissible kernels $k$,
their normalized Grammians are uniformly bounded below:
$$
N G^k \ \geq \  I
$$
for some $N$.

(iii) The sequence  $\{ \l_i \}_{i=1}^\i $ is strongly separated and
condition (a) alone holds.

(iv) Condition (b) alone holds.


Moreover, Condition (a) is equivalent to 

\noindent
(a'): There exists a constant $M$ and positive semi-definite infinite
matrices $\G$ and $\Delta$ such that
\be
\label{eq:bd}
M \delta_{ij} - 1 \= \G_{ij} (1 - \bar \l_i^1 \l_j^1) + \Delta_{ij}
(1 - \bar \l_i^2 \l_j^2) .
\ee

Condition (b) is equivalent to

\noindent
(b'): There exists a constant $N$ and positive
semi-definite infinite
matrices $\G'$ and $\Delta'$ such that
\be
\label{eq:bb}
N  -   \delta_{ij}  \= \G_{ij}' (1 - \bar \l_i^1 \l_j^1) + \Delta_{ij}'
(1 - \bar \l_i^2 \l_j^2) .
\ee

\et

We prove Theorem~\ref{thm:is} in Sections~\ref{sec:i-ii} and
\ref{sec:vec}.

\vskip 5pt

%

Compare Theorem~\ref{thm:is} with the recent characterization, by D.~Marshall and C.~Sundberg
and independently by C.~Bishop, of the interpolating sequences for the multiplier
algebra of the Dirichlet space, the space of analytic functions on
$\D$ with finite Dirichlet integral. The reproducing kernel for this space is
the Dirichlet kernel, $\displaystyle - \frac{1}{\bar \l z} \log (1 - \bar \l
z)$.

\bt {\bf [Marshall-Sundberg, Bishop] }

A sequence $\{ \l_i \}_{i=1}^\i $
in $\D$ is an interpolating sequence for the multiplier algebra of the
Dirichlet space if and only if it is 
weakly separated and the normalized Grammian of the Dirichlet kernel is
bounded.
Moreover, there are strongly separated sequences that are not interpolating
sequences.
\et

Possible strengthenings
of Theorem~\ref{thm:is} remain open, for example:

\begin{question}
If $\{ \lambda_i \}_{i=1}^\i $ is strongly separated, must it be an
interpolating sequence for $H^\i(\D^2)$?
\end{question}

\section{Background on kernels}
\label{sec:ker}
By a {\it kernel} $k$ on a set $X$, we mean a function
$k: X \times X \to \C$ that is positive definite, in the sense that
$$
\sum_{i,j = 1}^N \overline{c_i} c_j k(\l_i, \l_j) \ > \ 0
$$ for all $\l_1,\dots,\l_N$ in $X$ and all complex numbers
$c_1, \dots, c_N$, unless $c_1 = c_2 = \dots = c_N = 0$.

We shall use $k_\l$ to denote the function $k(\l,\cdot)$, and call this
function the kernel function at $\l$.

On the vector space of finite linear combinations of kernel
functions, \ie sums of the
form $\sum c_i k_{\l_i}$, one can define an inner product by defining
$$
\la k_\l , k_z \ra \= k(\l,z)
$$
and extending by sesqui-linearity.
Completing this vector space with this inner product gives a Hilbert
space of functions on $X$,
which we shall denote by $\hk$. For details of this construction,
see \cite{aro50}. Note that for any function $f$ in $\hk$
the construction yields the reproducing property of the kernel:
$$
\la f, k_\l \ra \= f(\l) .
$$

The multiplier algebra of $\hk$, denoted $M(\hk)$, is the set of
functions $\phi$ on $X$ with the property that whenever $f$ is in
$\hk$, then so is $\phi f$. It follows from the closed graph theorem
that if $\phi$ is a multiplier of $\hk$, then multiplication by $\phi$
is a bounded operator on $\hk$, and the operator norm makes $M(\hk)$
into a Banach algebra. We shall always consider $M(\hk)$ with this norm.

The most well-known non-trivial example is 
the \sz kernel. The Hilbert space produced is the Hardy space $H^2$, and
its multiplier algebra is (isometrically) $\hd$.

If $\phi$ is in $M(\hk)$, let $M_\phi$ denote the operaor on $\hk$ of
multiplication by $\phi$. Notice that all the kernel functions are
eigenvectors for the adjoint:
$$
M_\phi^* k_\l = \overline{\phi(\l)} k_\l ,
$$
as is seen by taking the inner product of either side with an arbitrary
function in $\hk$.

Notice too:
\begin{eqnarray*}
& \V M_\phi \V & \leq 1 \\
\Leftrightarrow& I - M_\phi M_\phi^* & \geq 0 \\
\Leftrightarrow&
(1 - \phi(z) \overline{\phi(\l)} ) k(\l,z) & \geq 0 .
\end{eqnarray*}

\vskip 5pt
Let us show that the boundedness of the Gram matrix
is equivalent to a Carleson measure condition. The result is well-known.

Let $k_i$ denote the kernel function at $\l_i$, and
$g_i$ denote the normalized kernel function at $\l_i$, so
$$
g_i \= \frac{1}{\V k_{\l_i} \V} k_{\l_i}  \= \frac{k_i}{\| k_i \|} .
$$

\bprop
\label{bgcm}
Let $\{ \l_i \}_{i=1}^\i $ be  a sequence in $X$. Then the following
conditions are equivalent.

(BG) The associated Gram matrix has norm at most $C$.

(CM) The measure $\sum \V k_i \V^{-2} \delta_{\l_i}$ is a Carleson measure
for $\hk$, \ie the
following inequality holds:
$$
\sum_{i=1}^\i \frac{|f (\l_i) |^2}{\V k_i \V^2} \ \leq \
C^2 \, \V f \V_{\hk}^2 .
$$
\eprop

\bp
(CM) $\Rightarrow$ (BG):
\begin{eqnarray*}
\V \sum a_i g_i \V_{\hk} \ &=& \
\sup_{\| f \| = 1} \la \sum a_i g_i , f \ra \\
&=& \sup_{\| f \| = 1} \sum a_i \| k_i \|^{-1} \overline{f(\l_i)}  \\
&\leq& \sup_{\| f \| = 1}
\left( \sum |a_i|^2 \right)^{1/2} \left( \sum \|k_i\|^{-2}
\ |f(\l_i)|^2 \right)^{1/2} \\
&\leq& C \, \left( \sum |a_i|^2 \right)^{1/2} .
\end{eqnarray*}

(BG) $\Rightarrow$ (CM):
Let $f$ be an arbitrary function in $\hk$, and
let $a_i = \| k_i \|^{-1} f(\l_i)$. Then
\begin{eqnarray*}
\sum \|k_i\|^{-2} |f(\l_i)|^2 \ &=& \
\la f , \sum a_i g_i \ra \\
&\leq& \| f \| \ \| \sum a_i g_i \| \\
&\leq& C \, \|f \| \  \left( \sum |a_i|^2 \right)^{1/2} .
\end{eqnarray*}
As $\sum |a_i|^2 = \sum \|k_i\|^{-2} |f(\l_i)|^2 $, we get the desired
inequality.
\ep

\section{Proofs of $(i) \Leftrightarrow (ii), (a) \Leftrightarrow (a')$ and
$(b) \Leftrightarrow (b')$}
\label{sec:i-ii}

First, let us prove that {\it (i)} and 
{\it (ii)} are equivalent in Theorem~\ref{thm:is}.

Given an interpolating sequence $\{ \l_i \}$, we shall call its
{\it interpolation constant} the infimum of those numbers $M$ such that, whenever
$|w_i| \leq 1$, there is a function $f$ of norm less than or equal to $M$
that interpolates each $\l_i$ to $w_i$.

\bl
The sequence $\{ \l_i \}_{i=1}^\i$ is an interpolating sequence for $\hb$
with interpolation constant $M$ if and only if, whenever $w_i$ is a sequence
of complex numbers with $\sup |w_i | \leq 1$ and $k$ is an admissible kernel,
then
\be
\label{eq:175}
M^2 \| \sum a_i k_{i} \|^2 \ \geq \
\| \sum a_i w_i k_i \|^2
\ee
whenever $\sum a_i k_{i}$ is in $\hk$.
\el
\bp
$(\Rightarrow)$.
By hypothesis, there is a function $\phi$ of norm less than or equal to $1$ in $\hb$
that maps $\l_i$ to $\bar w_i/M$. As $M_{z^1}$ and $M_{z^2}$ are commuting
contractions on $\hk$, by And\^o's inequality \cite{and63}, $\phi(M_{z^1},
M_{z^2}) = M_\phi$
is a contraction on $\hk$. Therefore
\begin{eqnarray*}
0 &\leq& M^2 \la \, (I - M_\phi  M_\phi^*)\, k_{\l_j} , k_{\l_i} \ra \\
&=& M^2  \la k_j , k_i \ra - w_j \bar w_i \la k_j , k_i \ra .
\end{eqnarray*}
But this means precisely that for any finite set of numbers $\{ a_i \}$,
we have
$$
M^2 \la \sum a_j k_j , \sum a_i k_i \ra \ \geq \ 
\la \sum a_j  w_j k_j , \sum a_i w_i k_i \ra ,
$$
and so (\ref{eq:175}) holds. 

$(\Leftarrow)$ 
Conversely, a necessary and sufficient condition 
to be able to find a function $\phi$ in $\hb$ of norm at most $M$ that
interpolates the points $\l_i$ to $\bar w_i$ is that, for every admissible
kernel $k$,
$$
(M^2 - w_j \bar w_i) \la k_j , k_i \ra \geq 0 .
$$
This was proved by the first author in \cite{ag1}; see also \cite{colwer94,
baltre98, agmc_bid}.
So if (\ref{eq:175}) holds for every admissible kernel and every choice of
$w_i$,
the sequence $\{ \l_i \}$ is
interpolating as desired.
\ep

Now, letting $w_j = \exp(2\pi i t_j)$  and $a_j = c_j$ in 
(\ref{eq:175}) and integrating with respect to $t_1, t_2, \dots$ on
$[0,1] \times [0,1] \times \dots$ one gets
\be
\label{eq:206}
\sum_j |c_j|^2 \| k_j \|^2 \ \leq \
M^2 \| \sum_j c_j k_j \|^2 .
\ee
Similarly, letting
$a_j = \exp(- 2\pi i t_j) c_j$ and $w_j = \exp(2\pi i t_j)$
and integrating gives
\be
\label{eq:208}
\| \sum_j c_j k_j \|^2 \ \leq \
M^2 \sum_j |c_j|^2 \| k_j \|^2 .
\ee
Combining (\ref{eq:206}) and (\ref{eq:208}), one gets that if $\{ \l_i \}$ is
an interpolating sequence, then for any normalized admissible kernel we have
\be
\label{eq:209}
\frac{1}{M^2} \sum_i |c_i|^2  \leq \| \sum_i c_i g_i \|^2 \leq M^2 \sum_i
|c_i|^2 
\ee
(in other words, $\{ g_i \}$ is a {\it Riesz system}, and the constant $M$
can be chosen uniformly).
Conversely, if (\ref{eq:209}) holds, then 
(\ref{eq:175}) holds, with the constant $M^2$ replaced by $M^4$. 
As the first inequality in (\ref{eq:209}) says $G$ is bounded below by
$1/M$, and the second inequality says $G$ is bounded by $M$, we have shown:

{\it The sequence $\{ \l_i \}_{i=1}^\i$ is an interpolating sequence with
interpolation constant $M$.

$\Rightarrow$
$\displaystyle \frac{1}{M} I \leq G^k \leq M I$ for all admissible kernels
$k$.

$\Rightarrow$
The sequence $\{ \l_i \}_{i=1}^\i$ is an interpolating sequence with
interpolation constant $M^2$.
}

\vskip 10pt
To show that $(a') \Rightarrow (a)$ and 
$(b') \Rightarrow (b)$ is easy: take the Schur product (the entrywise
product) of both sides of (\ref{eq:bd}) with $G^k$, and one gets:
$$
M I - G^k = \Gamma \cdot [ (1 - \bar \l_i^1 \l_j^1) \cdot G^k ]
+ \Delta \cdot [ (1 - \bar \l_i^2 \l_j^2) \cdot G^k ] .
$$
As the two quantities in brackets are positive, by definition of an
admissible kernel, and the Schur product of two positive matrices is
positive, one gets $ G^k \leq MI$.
Similarly, (\ref{eq:bb}) gives $NG^k \geq I$.

The  converse direction is a duality argument.
Suppose $(a)$ holds. This can be expressed as saying:

Whenever $k(\l,z)$ is a kernel such that:
\setcounter{equation}{\value{theorem}}
\begin{eqnarray}
(1 - \overline{\l^1} z^1) \cdot k \ &\geq& \ 0 
\label{eq:2451} \\
\quad\mbox{and}  \qquad
(1 - \overline{\l^2} z^2) \cdot k \ &\geq& \ 0 
\label{eq:2452} \\
\mbox{then} && \nonumber \\
(MI - J) \cdot k \ &\geq& \ 0  \nonumber ,
\addtocounter{theorem}{2}
\end{eqnarray}
where we use $\cdot$ to denote the Schur product.
Now a Hahn-Banach argument shows that $MI - J$ cannot be separated from the
closed wedge of infinite matrices of the form
$$
\Gamma(\l,z) \cdot  (1 - \bar \l^1 z^1) + \Delta(\l,z) \cdot  (1 - \bar \l^2 z^2)
\quad : \quad \Gamma \geq 0,\, \Delta \geq 0 .
$$

Indeed, fix a positive integer $n$. 
Let ${\cal T}_n$ be the set of all $n$-by-$n$ self-adjoint matrices  $T$
representable
in the form
$$
T_{ij} \= (1 - \lambda^1_i \overline{\lambda^1_j} ) \Gamma_{ij} +
(1 - \lambda^2_i \overline{\lambda^2_j} ) \Delta_{ij} ,
$$
where $\Gamma$ and $\Delta$ are positive.
As ${\cal T}_n$ is a closed wedge, if $MI - J$ were not in ${\cal
T}_n$, there
would
be a real linear functional on the space of all $n$-by-$n$
self-adjoint matrices
that was positive on ${\cal T}_n$ and strictly negative on $MI - J$.

Any such linear functional must be of the form
$T \mapsto tr(WT)$ for some self-adjoint matrix $W$.
As ${\cal T}_n$ contains the set of all positive matrices (let $\Delta =
0$ and
$\Gamma$ be the Schur product of the positive matrix
$1/(1 - \overline{\lambda^1_i} \lambda^1_j)$
with an arbitrary positive matrix), $W$ must be positive.
Let $K$ be the transpose of $W$. Then
$$
\sum_{i,j=1}^n c_i \overline{c_j} K_{ij} (1 - \l^1_i
\overline{\l^1_j})
\= tr( WT) ,
$$
where
$$
T_{ij} \= (1 - \l^1_i \overline{\l^1_j} ) c_i \overline{c_j}
$$
is in ${\cal T}_n$. So $K$ satisfies (\ref{eq:2451}), and similarly also
(\ref{eq:2452}).
Therefore $(MI - J) \cdot K$ is
positive, so
$$
tr(W (MI -J)) \= \sum_{i,j=1}^n K(\lambda_i, \lambda_j) (M
\delta_{ij} -1) \ \geq
0 ,
$$
a contradiction if $MI - J$ is not in ${\cal T}_n$.
So for every $n$, we have $MI - J$ is in ${\cal T}_n$;
it follows \eg from Kurosh's theorem \cite[p.75]{arkh} that there is
a choice of $\Gamma$ and $\Delta$ such that
$$
M \delta_{ij} - 1 \= \G_{ij} (1 - \bar \l_i^1 \l_j^1) + \Delta_{ij}
(1 - \bar \l_i^2 \l_j^2)
$$
for all $i,j$.

A similar argument shows that $(b) \Rightarrow (b')$.

\section{Proof that $(iii) \Leftrightarrow (iv)$.}
\label{sec:vec}

Let us analyze condition $(b')$.
First some notation. Given Hilbert spaces
$\L_1$ and $\L_2$, we let $B(\L_1,\L_2)$ denote the bounded 
linear operators from $\L_1$ to $\L_2$, and $H^\i(\D^2,
B(\L_1,\L_2))$ the space of bounded holomorphic functions from $\D^2$
to $B(\L_1,\L_2)$. Let $\{e_i \}_{i=1}^\i$ be the usual orthonormal
basis of $l^2$, the vector with $1$ in the $i^{th}$ slot and $0$
elsewhere.

\bl
\label{lem:l2}
With notation as in Theorem~\ref{thm:is}, condition (b') is
equivalent to:

(b''): There exists a function $\Phi$ in $H^\i(\D^2,B(\C,l^2))$ of
norm at most $\sqrt{N}$ such that 
$ \Phi(\l_i) = e_i$.
\el

Before proving this lemma, we need to recall the following theorem of
the first author \cite{ag90}:

\bt
\label{thm:rep}
 The function $\Psi$ is in the closed unit ball 
of $H^\i(\D^2, B(\L_1,\L_2))$ if and only if: there are auxiliary
Hilbert spaces $\h_1$ and $\h_2$ and an isometry
$U:\L_1 \oplus \h_1 \oplus \h_2 : \rightarrow \L_2 \oplus \h_1 \oplus
\h_2 $ such that, with respect to the decomposition of $U$ as
$$
U \= \bordermatrix{&\L_1 &\h_1 \oplus \h_2\cr
\L_2 &A & B\cr
\h_1 \oplus \h_2&C &D}
$$
we have
$$
\Psi(\lambda) \= A + B E_\lambda (1 - D E_\lambda)^{-1} C .
$$
Here, for $\lambda = (\l^1,\l^2)$ in $\D^2$,
$\displaystyle E_\lambda \, = \, \lambda^1 I_{\h_1} \oplus
 \lambda^2 I_{\h_2} $ is the operator of multiplication by $\l^1 $ on $\h_1$
and
multiplication by $\l^2$ on $\h_2$.
\et

\vskip 5pt
{\sc Proof of Lemma~\ref{lem:l2}:}
Consider condition $(b')$:
\be
\label{eq:4b}
N  -   \delta_{ij}  \= \G_{ij}' (1 - \bar \l_i^1 \l_j^1) +
\Delta_{ij}'
(1 - \bar \l_i^2 \l_j^2) .
\ee
Choose vectors $f_i$ and $g_i$ in auxiliary Hilbert spaces $\h_1$
and $\h_2$ so that
\begin{eqnarray*}
\la f_j, f_i \ra &\=& \G_{ij}' \\
\la g_j, g_i \ra &\=& \Delta_{ij}' .
\end{eqnarray*}
Then Equation~\ref{eq:4b} can be rewritten as
\be
\label{eq:2651}
N + \bar \l_i^1 \l_j^1 \la f_j, f_i \ra + \bar \l_i^2 \l_j^2 \la g_j, g_i
\ra \= \la e_j, e_i \ra + \la f_j, f_i \ra + \la g_j, g_i \ra.
\ee
Letting $h_i = f_i \oplus g_i$, Equation~\ref{eq:2651} becomes
\be
\label{eq:2652}
\la \left(\matrix{\sqrt{N}\cr E_{\lambda_j} h_j\cr } \right) ,
\left(\matrix{\sqrt{N}\cr E_{\lambda_i} h_i\cr } \right) \ra
\=
\la \left(\matrix{e_j\cr h_j\cr } \right) ,
\left(\matrix{e_i\cr h_i\cr } \right) \ra .
\ee
So there is an isometry
$$
L: \left(\matrix{\sqrt{N}\cr E_{\lambda_i} h_i\cr } \right) \mapsto
\left(\matrix{e_i\cr h_i\cr } \right).
$$
Increasing $\h_1$ and $\h_2$  if necessary, $L$ can be extended to a unitary
$U: \C \oplus \h_1 \oplus \h_2 \rightarrow l^2 \oplus \h_1 \oplus \h_2$.
Write 
$$
U \= \bordermatrix{&\C &\h_1 \oplus \h_2\cr
l^2 &A & B\cr
\h_1 \oplus \h_2&C &D} ,
$$
and let $\Psi(\lambda) \= A + B E_\lambda (1 - D E_\lambda)^{-1} C$.
Then, by Theorem~\ref{thm:rep}, $\V \Psi \V \leq 1$.
Solving
\be
\label{1118}
\left(\matrix{A & B \cr
C & D \cr} \right)  \left(\matrix{\sqrt{N}\cr E_{\lambda_i} h_i\cr }
\right) \=
\left(\matrix{e_i\cr h_i\cr } \right),
\ee
we get $\Psi(\l_i)
= \frac{1}{\sqrt{N}} e_i$. Then $\Phi = \sqrt{M} \Psi$ is the required
function.

Conversely, if $(b'')$ holds, let $\Psi = \frac{1}{\sqrt{N}} \Phi$,
and write $\Psi$ as in Theorem~\ref{thm:rep}. Then
Equation~\ref{1118} holds, and hence going backwards so do
Equations~\ref{eq:2652}, \ref{eq:2651} and \ref{eq:4b}.
\ep

A similar argument shows:
\bl
\label{lem:l3}
With notation as in Theorem~\ref{thm:is}, condition (a') is
equivalent to:

(a''): There exists a function $\Psi$ in $H^\i(\D^2,B(l^2,\C))$ of
norm at most $\sqrt{M}$ such that
$ \Psi(\l_i)  e_i = 1 $.
\el

\vskip 5pt

Now, suppose condition $(b'')$ holds.
Letting $\Psi(\lambda) = \Phi(\lambda)^t$, we get $(a'')$.
Moreover, writing $\Phi$ as
$$
\Phi(\lambda)  \= \left(\matrix{\phi_1 \cr \phi_2 \cr \vdots \cr}
\right),
$$
we get a sequence of functions $\phi_i$ such that
$\displaystyle \sum_i |\phi_i (\lambda)|^2 \leq N$
and $\phi_i(\l_j) = \delta_{ij}$. So in particular,
$\{ \l_i \}$ is strongly separated, proving that $(iv) \Rightarrow (iii)$.

Conversely, suppose $(iii)$ holds. By $(a'')$, writing
$$
\Psi = (\psi_1,\psi_2,\dots) ,
$$
we have $\displaystyle \sum_i |\psi_i (\lambda)|^2 \leq M$
and $\psi_i(\l_i) = 1$.
Moreover, strong separation means we have a sequence of functions $\chi_i$
such that $\chi_i(\l_j) = \delta_{ij}$ and $\| \chi_i \| \leq C$ for all
$i$.
Letting $\phi_i = \psi_i \chi_i$,
and 
$$
\Phi(\lambda)  \= \left(\matrix{\phi_1 \cr \phi_2 \cr \vdots \cr}
\right),
$$
we get $\V \Phi \V \leq C \sqrt{M}$, and $\Phi(\l_i) = e_i$.
So $(b'')$ holds with constant $C \sqrt{M}$, proving $(iii) \Rightarrow
(iv)$.
\ep
\vskip 5pt
Remark: 
A theorem due to Varopoulos \cite{var71}
and Bernard \cite{ber71} (see also \cite[p. 298]{gar81}), which applies
to any
uniform algebra that is also a dual space, asserts that given an
interpolating sequence $\{ \l_i \}_{i=1}^\i$, one can find functions
$\phi_i$
with $\phi_i(\l_j) =
\delta_{ij}$ and $\sum |\phi_i (z) | \leq M^2$ for all $z$. These are
sometimes
called Per Beurling functions, because he showed they existed for all
interpolating sequence on the disk (see \cite[p. 294]{gar81}).
The equivalence of Condition~$(b'')$ and a sequence being
interpolating is then the $\hb$ case of 
the Varopoulos-Bernard theorem.

\end{document}